\documentclass[a4paper, 12pt]{amsart}
\usepackage[dvips]{graphicx}
\usepackage[all]{xy}

\newtheorem{thm}{Theorem}[section]
\newtheorem{prop}[thm]{Proposition}
\newtheorem{lem}[thm]{Lemma}
\newtheorem{definition}[thm]{Definition}
\newtheorem{corollary}[thm]{Corollary}

\newtheorem{theorem}{Theorem}[section]

\newtheorem{remark}{\bf Remarks}[section]

\numberwithin{equation}{section} 

\newcommand{\bea} {\begin{eqnarray*}}
\newcommand{\beq} {\begin{equation}}
\newcommand{\bey} {\begin{eqnarray}}
\newcommand{\eea} {\end{eqnarray*}}
\newcommand{\eeq} {\end{equation}}
\newcommand{\eey} {\end{eqnarray}}

\begin{document}
\title[Minimal fixed point set of maps on T-Bundles over $S^1$]
{Minimal fixed point set of maps on Torus Fiber Bundles over the Circle}
\author{Weslem L. Silva}
\address[Weslem Liberato Silva]{Deptartment of Mathematics and Computer Science, Loyola University, 
6363 St charles Avenue, New Orleans, LA 70118, U.S.A.}
\email{weslemliberato@gmail.com}
\begin{abstract}
The main purpose this work is to study the minimal fixed point set 
of fiber-preserving maps for spaces which are fiber  bundles  over the circle and the fiber is the  torus. 
Using the one-parameter fixed point theory is possible to describe these sets in terms of the fundamental group and the induced 
homomorphism.
\end{abstract}
\date{\today}
\maketitle
\bibliographystyle{amsplain}

\footnotetext[1]{This work is part of the doctoral thesis of the author, made at University Federal of S\~ao
Carlos under the supervision of Professor Daniel Vendr\'uscolo. 

This work was supported by CAPES.}


\section{Introduction}
Let $S \to M \stackrel{p}{\to} B$ be a fiber bundle, where $S,M,B$ are closed manifolds, and $f: M \to M$ be a 
fiber-preserving map . The minimum number $MF_{B}[f] = min \{ \# \pi_{0}(Fix(f^{'})) | f^{'} \sim_{B} f \}$ 
of path components of fixed point subspaces of $M$ among all pairs fiberwise homotopic to $f$ is finite, see \cite{K-11}.
The symbol \lq\lq $\sim_{B}$ \rq\rq \, means a fiberwise homotopy.

To determine when the number $MF_{B}[f]$ is zero, that is, when the fiber-preserving map 
$f$ can be deformed by a fiberwise homotopy to a fixed point free map is a problem that has been
considered by many authors, see for example, \cite{F-H-81}, \cite{G-P-V-04} and \cite{S-V-12}.
The study of the minimal fixed point set of a fiber-preserving map is a problem of interest  
in fixed point theory. These sets have been studied using bordism techniques,  
that in general are difficult to compute, see \cite{K-11}.

In this paper we present a method 
to compute $MF_{B}[f]$, using one-parameter fixed point theory, when the base $B$ is the circle $S^{1}$. This technique allows us to present 
the minimal fixed point set of a fiber-preserving map in terms of the fundamental group of $M$, and 
of the induced homomorphism $f_{\#}$. The one-parameter fixed point theory also allow us to describe each path component 
of $Fix(f^{'})$ for each fiber-preserving map $f^{'}$ fiberwise homotopic to $f$.

Let $f : M \to M$ be a fiber-preserving map, where $M$ is a fiber bundle
over the circle  and the fiber is the torus, $T$. Such fiber bundles $M$ are
obtained from $T \times [0,1]$ by identifying $(x, 0)$ with $(A(x),1)$, where $A$ is a
homeomorphism of T. We write $$ M(A) = M = \frac{T \times [0,1]}{(x,0) \sim (A(x),1)} $$
The elements of $M(A)$ are denotes by $<[(x,y)],t>$. Here $[(x,y)]$ denote a point in $T$. 
We can identify $A$ with a matrix with integer coefficients and determinant 
$1$ or $-1$, the details are in section \ref{section-3}. The projection map 
$p: M(A) \to S^{1} = I/0 \sim 1 $, is given by $p(<[(x,y)],t>)$ $ = <t>$.

Since $f$ is a fiber-preserving map and the base is $S^{1}$, the fixed point set of $f$ can be seen as 
the fixed point set of a homotopy of the torus.  In this paper we study the minimal fixed point set for 
homotopies using one-parameter fixed point theory developed by  R. Geoghegan and A. Nicas in 
\cite{G-N-94}.

This paper is organized into five sections, besides this one. In section 2 we considered fiber-preserving maps in fiber 
bundles over the circle with fiber torus. In section 3 we present the relation between fixed point sets of fiber-preserving maps 
and fixed point sets of homotopies. In section 4 we present  preliminares 
about the one-parameter fixed point theory. 
In section 5 we prove the main result, which is theorem 
\ref{maintheorem}.


\section{Torus fiber-preserving maps} \label{section-3}

Let $T$ be, the torus, defined as the quotient space $ {\mathbb{R} \times \mathbb{R}}/{\mathbb{Z} \times \mathbb{Z}} $. 
We denote by $(x, y)$ the elements of  $\mathbb{R} \times \mathbb{R}$ and by  
$[(x,y)]$ the elements in T.

Let $M(A) = \frac{T \times [0,1]}{([(x,y)],0) \sim (\left[A \left(^{x}_{y} \right) \right],1)}  $
be the quotient space, where $A$ is a homeomorphism of $T$ induced by an operator in 
$\mathbb{R}^{2}$ that preserves $ \mathbb{Z} \times \mathbb{Z}$. 
The space $M(A)$ is a fiber bundle over the circle $S^{1}$ where the fiber is the torus.
For more details on these bundles see \cite{G-P-V-04}.

Given a fiber-preserving map $f: M(A) \to M(A) $, i.e. $p \circ f = p$ we want to compute the number $MF_{S^{1}}[f]$. 
More precisely we want to study the path components of $Fix(f^{'})$ for each map $f^{'}$ fiberwise homotopic to $f$.

Consider the loops in $M(A)$  given by;
$a(t) = <[(t,0)],0> $, $b(t) = $ $ <[(0,t)],0> $ and 
$c(t) = <[(0,0)],t> $ for $t \in [0,1]$.
We denote by $B$ the matrix of the homomorphism induced on the fundamental group by the 
restriction of $f$ to the fiber $T.$ From \cite{G-P-V-04} we have the following theorem
that provides a relationship between the matrices $A$ and $B$, where 
$$ A = { \left( \begin{array}{cc} a_{1} & a_{3} \\ a_{2} & a_{4} \\ \end{array} \right)}  $$
From \cite{G-P-V-04} the induced homomorphism $f_{\#}: \pi_{1}(M(A)) \to \pi_{1}(M(A)) $ is given by 
$f_{\#}(a) = a^{b_{1}} b^{b_{2}} $, $f_{\#}(b) = a^{b_{3}} b^{b_{4}} $, $f_{\#}(c) = a^{c_{1}} b^{c_{2}} c $.  
Thus $$ B = { \left( \begin{array}{cc} b_{1} & b_{3} \\ b_{2} & b_{4} \\ \end{array} \right)} $$

\begin{theorem} \label{daci1theorem}
$(1) \,\, \pi_{1}(M(A),0) = \langle a,b,c | [a,b] = 1, cac^{-1} = a^{a_{1}} b^{a_{2}} , 
 cbc^{-1}  $ $ = a^{a_{3}} b^{a_{4}} \rangle $ 

$(2)\, B $ commutes with $A$.

$(3) \, $ If $f$ restricted to the fiber is deformable to a fixed point free map then 
the determinant of $B - I$ is zero, where $I$ is the identity matrix.

$(4) \, $ Consider $w = A(v)$ if the pair $v,w$ generators $ \mathbb{Z} \times \mathbb{Z}$, 
otherwise let $w$ be another vector so that $v,w$ span $\mathbb{Z} \times \mathbb{Z}$.
Define the linear operator 
$P:\mathbb{R} \times \mathbb{R} \to \mathbb{R} \times \mathbb{R} $ by $P(v) = \left( ^{1}_{0} \right)$
and  $P(w) = \left( ^{0}_{1} \right)$. Consider an isomorphism of fiber bundles, also denoted by $P$,
$P: M(A) \to M(A^{1})$ where $A^{1} = P\circ A \circ P^{-1}$. Then $M(A)$ is homeomorphic to $M(A^{1})$ 
over $S^{1}$. Moreover we have one of the cases of the table below with $B^{1} = P \circ A \circ P^{-1}$
and $B \neq I$, except in case $I$:

{\small \begin{center}
\begin{tabular}{|p{1,3cm}|l|}
\hline
 $ Case \,\, I $ & $ A^{1} = \left( \begin{array}{cc} a_{1} & a_{3} \\ a_{2} & a_{4} \\ \end{array} \right) $ ,
 $ B^{1} = \left( \begin{array}{cc} 1 & 0 \\ 0 & 1 \\ \end{array} \right) $ \\
  & $a_{3} \neq 0 $    \\
\hline
$ Case \,\, II $ & $ A^{1} = \left( \begin{array}{cc} 1 & a_{3}\\ 0 & 1 \\ \end{array} \right) $ ,
 $ B^{1} = \left( \begin{array}{cc} 1 & b_{3} \\ 0 & b_{4} \\ \end{array} \right) $ \\
  & $a_{3}(b_{4}-1) = 0 $    \\
\hline
$ Case \,\, III $ & $ A^{1} = \left( \begin{array}{cc} 1 & a_{3} \\ 0 & -1 \\ \end{array} \right) $ ,
 $ B^{1} = \left( \begin{array}{cc} 1 & b_{3} \\ 0 & b_{4} \\ \end{array} \right) $ \\
  & $a_{3}(b_{4}-1) = -2b_{3} $    \\
\hline
$ Case \,\, IV $ & $ A^{1} = \left( \begin{array}{cc} -1 & a_{3} \\ 0 & -1 \\ \end{array} \right) $ ,
 $ B^{1} = \left( \begin{array}{cc} 1 & b_{3} \\ 0 & b_{4} \\ \end{array} \right) $ \\
  & $a_{3}(b_{4}-1) = 0 $    \\  
\hline
$ Case \,\, V $ & $ A^{1} = \left( \begin{array}{cc} -1 & a_{3} \\ 0 & 1 \\ \end{array} \right) $ ,
 $ B^{1} = \left( \begin{array}{cc} 1 & b_{3} \\ 0 & b_{4} \\ \end{array} \right) $ \\
  & $a_{3}(b_{4}-1) = 2b_{3} $    \\  
\hline
\end{tabular} \end{center}} 
\end{theorem}

From \cite{G-P-V-04} we have the following theorem:

\begin{theorem} \label{main-theorem-daci1}
If $f: M(A) \to M(A) $ is a fiber-preserving map, then in the case $I$ we have $MF_{S^{1}}[f]= 0$,
and in the cases $II$ and $III$ we have $MF_{S^{1}}[f]= 0$  if and only if $c_{1}(b_{4}-1)-c_{2}b_{3}=0$.

\end{theorem}

The Theorem \ref{main-theorem-daci1} in \cite{G-P-V-04} provides also conditions for remaining cases. We omit them, since here we will study only $II$ and $III$.



\section{Fixed point set of fiber-preserving maps} \label{section-1}

Given a fiber-preserving map $f: M(A) \to M(A) $  the set 
$Fix(f)$ is given by; $\{ \,\, <[(x,y)],t> \, \in M(A) \,\, | \,\, f(<[(x,y)],t>) = <[(x,y)],t> \,\, \}$.
Since $f$ is a fiber-preserving map then 
the map $f$ is given by formula: $$f(<[(x,y)],t>) = <F([(x,y)],t),t >  $$
where $F: T \times I \to T$ is a homotopy. 
We call this homotopy $F$ the homotopy induced by $f$. 
If $f$ has no fixed points in $t=0,1,$ then the study of the set $Fix(f)$ is equivalent to the study of the set $Fix(F)$, that is, 
$$Fix(f) \approx Fix(F). $$
This happens since, in the fiber bundle $M(A)$ the class $<[(x,y)],t>$ contains only one unique point if $t \neq 0,1$.  
Notice that

\begin{prop} \label{defproposition}
Let $M(A)$ be a fiber bundle as in theorem \ref{daci1theorem}. If $f: M(A) \to M(A)$ is 
a fiber-preserving map such the restrction to each fiber $f_{|_{T}}$ can be deformed to a fixed point free map,
then $f$ can be deformed to a map $f^{'}$ such that 
$f^{'}(<[(x,y)],0>): T \to T$ is a fixed point free map.
\end{prop}
{\it Proof.}
Let $f: M(A) \to M(A)$ be a fiber-preserving map given by $f(<[(x,y)],t>) = <F([(x,y)],t),t>$.
As $M(A)$ is a locally trivial bundle thus we can choose $\frac{1}{2} > \epsilon > 0$ 
such that $p^{-1}((\epsilon, 1-\epsilon)) \approx   T \times (\epsilon, 1-\epsilon) $. 
We take the homotopy
$H: M(A) \times I \to M(A)$ defined by;
$$ { H(<[(x,y)],t>,s) =  \left \{  
\begin{array}{lll}
 <F([(x,y)],0),t> & if & 0 \leq t \leq s\epsilon \\
 <F([(x,y)],\frac{1}{1-2s\epsilon}(t-s\epsilon)),t> & if & s\epsilon \leq t \leq 1-s\epsilon \\
 <F([(x,y)],0),t> & if & 1-s\epsilon \leq t \leq 1 \\
\end{array} \right.}
$$

By hypothesis there is one homotopy 
$h:T \times I \to T$ satisfying $h([(x,y)],1) = F([(x,y)],0)$ and  
$h([(x,y)],0)$ is a fixed point free map.
Therefore we can define the following homotopy:
$$ { \footnotesize G(<[(x,y)],t>,s) =  \left \{  
\begin{array}{lll}
 <h([(x,y)],\frac{(t-\epsilon)}{\epsilon}s+1),t> & if & 0 \leq t \leq \epsilon \\
 <F([(x,y)],\frac{1}{1-2\epsilon}(t-\epsilon)),t> & if & \epsilon \leq t \leq 1-\epsilon \\
 <h([(x,y)],\frac{-(t-1+\epsilon)}{\epsilon}s+1),t> & if & 1-\epsilon \leq t \leq 1 \\
\end{array} \right. }
$$
The fiber-preserving homotopy $J:M(A) \times I \to M(A)$ defined by; 
$$ J(<[(x,y)],t>,s) =  \left \{  
\begin{array}{lll}
 H(<[(x,y)],t>,2s) & if & 0 \leq s \leq \frac{1}{2} \\
 G(<[(x,y)],t>,2s-1)& if & \frac{1}{2} \leq s \leq 1 \\
\end{array} \right.
$$ 
satisfies the condition of the theorem. \qed

\bigskip

Note that if a fiber-preserving map $f: M(A) \to M(A)$, where $M(A)$ is as in Theorem \ref{daci1theorem}, has no fixed points in $t = 0 $ 
then $f$ has no fixed points in $t = 1$ also. In fact, suppose that $f$ has one fixed point in $t=1$. We have,
$f(<[(x,y)],0> ) = f(<[A(^{x}_{y})],1>)$. Using which the matrix $A$ is invertible in $\mathbb{Z}$, see section \ref{section-3}, 
then there should be one  point $<[(u,v)],0>$, satisfying $f(<[(u,v)],0>) = <[(u,v)],0>$, but this is a contradiction.

\begin{prop} \label{propconjugation}
Let $F: T \times I \to T$ be the homotopy induced by a fiber-preserving map $f: M(A) \to M(A)$, i.e, 
$f(<[(x,y)],t>) = $  $<F([(x,y)],t),t>$. If $P: T \to T $  is an isomorphism and  
$g: M(A^{1}) \to M(A^{1})$, $A^{1} = P \circ A \circ P^{-1}$, is a fiber-preserving map  
defined by  $g(<[(x,y)],t>) = $  $ <P \circ F \circ (P^{-1} \times Id)\left(^{x}_{y},t 
\right), t>$, then the numbers $MF_{S^{1}}[f] $ and $ MF_{S^{1}}[g]$ are equals.
$$
\xymatrix{ M(A)  \ar[r]^{f} \ar[d]_{\bar{P}} &   M(A) \ar[d]^{\bar{P}}  \\
M(A^{1}) \ar[r]^{g}  &   M(A^{1}) \\ }
$$
\end{prop}
{\it Proof.} 
Note that the homotopy $G: T \times I \to T$ induces the fiber-preserving map $g: M(A^{1}) \to M(A^{1})$. 
Since that $G = P \circ F \circ (P^{-1} \times Id) $ then we have $MF_{S^{1}}[f] = MF_{S^{1}}[g]$. \qed

\bigskip

\begin{figure}[!htp]
\begin{center}
		\includegraphics[scale=0.28]{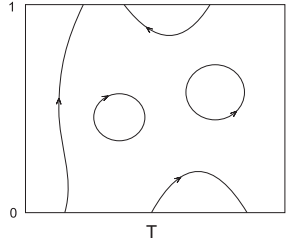}
	\caption{Fixed point set of a homotopy on the torus.}
	\label{fig:Figura13}
	\end{center}
\end{figure}

By Proposition \ref{defproposition}, the study of the minimal fixed point set, over $S^{1}$, of a fiber-preserving map $f: M(A) \to M(A)$ is 
equivalent to study of the minimal fixed point set for the homotopy induced by $f$. 
In this paper, we applied the one-parameter fixed point theory developed by R. Geoghegan and A. Nicas in 
\cite{G-N-94} to determine the minimal fixed point set of the homotopy $F$.
Since $T$ is orientable then the fixed point set of $F$ consists of 
oriented arcs as Figure \ref{fig:Figura13}, see \cite{D-G-90}, \cite{G-N-S-00} and \cite{Sc-83}.

As $f$ has no fixed points in $t=0,1,$ then in Figure \ref{fig:Figura13} above we have only circles. The one-parameter trace 
give us the \lq\lq minimum amount \rq\rq \, these circles. In the next section we shall describe some concepts of the one-parameter 
fixed point theory, for more details see \cite{G-N-94}.


\section{One-parameter fixed point theory} \label{section-2}

\subsection{Hochschild Homology Traces} 

Let $R$ be a ring and $M$ an $R-R$ bimodule, that is, a left and right R-module satisfying $(r_{1}m)r_{2} = r_{1}(mr_{2})$ 
for all $m \in M$, and $r_{1},r_{2} \in R$. The Hochschild chain complex $\{C_{\ast}(R, M) ,d \}$ is given by 
$C_{n}(R,M) = R^{\otimes n} \otimes M$ where $R^{\otimes n}$ is the tensor product of n copies of $R$, taken 
over the intergers, and  
$$
\begin{array}{lcl}
d_{n}(r_{1} \otimes \ldots \otimes r_{n} \otimes m ) & = &
 r_{2} \otimes \ldots r_{n} \otimes m r_{1}  \\
 &  & + \displaystyle \sum_{i=1}^{n-1} (-1)^{i} r_{1} \otimes \ldots \otimes 
 r_{i} r_{i+1} \otimes \ldots \otimes r_{n} \otimes m  \\
 &  & + \, (-1)^{n} r_{1} \otimes \ldots \otimes r_{n-1} \otimes r_{n} m . \\
\end{array}
$$

The $n-th$ homology of this complex is the Hochschild homology of $R$ with coefficient bimodule $M$. It 
is denoted by $HH_{n}(R,M)$. For computed $HH_{1}$ and $HH_{0}$ we have the formulae     
$ d_{2}(r_{1} \otimes r_{2}\otimes m) = r_{2}\otimes m r_{1} -  r_{1}r_{2}\otimes m + r_{1}\otimes r_{2}m$ and  
$ d_{1}(r \otimes m) = mr - rm $. 

\begin{lem}
If $1 \in R$ is the unit element and $m \in M$ then the 1-chain $1 \otimes m$ is a boundary.
\end{lem}
{\it Proof.} $d_{2}(1\otimes 1 \otimes m) = 1\otimes m - 1\otimes m + 1\otimes m = 1 \otimes m$. \qed

\bigskip

The Hochschild homology will arise in the following situaton: let $G$ be a group and $\phi:G \to G$ an endomorphism. 
Also denote by $\phi$ the induced ring homomorphism $\mathbb{Z}G \to \mathbb{Z}G$. 
Take the ring $R = \mathbb{Z}G$ and $M = (\mathbb{Z}G)^{\phi}$ the $\mathbb{Z}G-\mathbb{Z}G$ bimodule whose underlying 
abelian group is $\mathbb{Z}G$ and the bimodule structure is given by $g.m = gm$ and  $m.g = m \phi(g)$.

Two elements $g_{1},g_{2}$ in $G$ are semiconjugate if and only if there exists $g \in G$ such that $g_{1} = g g_{2} \phi(g^{-1})$.
We write $C(g)$ for the semiconjugacy class containing $g$ and $G_{\phi}$ for the set of semiconjugacy classes. 
Thus, we can decompose $G$ in the union of its semiconjugacy classes. 
This partition induces a direct sum decomposition of $HH_{\ast}(\mathbb{Z}G,(\mathbb{Z}G)^{\phi})$.

In fact, each generating chain $\gamma = g_{1} \otimes ... \otimes g_{n} \otimes m $ can be written in canonical form 
as $g_{1} \otimes ... \otimes g_{n} \otimes g_{n}^{-1} \otimes ... g_{1}^{-1}g $ where $g = g_{1} ... g_{n}m \in G$ ``marks'' 
a semiconjugacy class. Thus, the decomposition $(\mathbb{Z}G)^{\phi} \cong \bigoplus_{C \in G_{\phi}} \mathbb{Z}C$ as a direct 
sum of abelian groups determines a decomposition of chains complexes $C_{\ast}(\mathbb{Z}G,(\mathbb{Z}G)^{\phi}) \cong 
\bigoplus_{C \in G_{\phi}} {C_{\ast}(\mathbb{Z}G,(\mathbb{Z}G)^{\phi})}_{C} $ where ${C_{\ast}(\mathbb{Z}G,(\mathbb{Z}G)^{\phi})}_{C} $ 
is the subgroup of $C_{\ast}(\mathbb{Z}G,(\mathbb{Z}G)^{\phi})$ generated by those generating chains whose markers lie in $C$. 
Thus we have the following isomorphism: $HH_{\ast}(\mathbb{Z}G,(\mathbb{Z}G)^{\phi}) \cong 
\bigoplus_{C \in G_{\phi}} {HH_{\ast}(\mathbb{Z}G,(\mathbb{Z}G)^{\phi})}_{C} $ where the summand ${HH_{\ast}(\mathbb{Z}G,(\mathbb{Z}G)^{\phi})}_{C} $ 
corresponds to the homology classes marked by the elements of $C$. This summand is called the $C-$component. 

Let $Z(h) = \{g \in G| h = gh\phi(g^{-1})  \}$ be the semicentralizer of $h \in G$. Choosing representatives $g_{C} \in C$, then we have 
the following proposition whose proofs is in \cite{G-N-94}:
\begin{prop}
Choosing representatives $g_{C} \in C$ then we have 
$$HH_{\ast}(\mathbb{Z}G,(\mathbb{Z}G)^{\phi}) \cong 
\bigoplus_{C \in G_{\phi}} {H_{\ast}(Z(g_{C}))}_{C} $$
where ${H_{\ast}(Z(g_{C}))}_{C} $ corresponds to the summand ${HH_{\ast}(\mathbb{Z}G,(\mathbb{Z}G)^{\phi})}_{C}$. 
\end{prop}  

Given a $m \times n $ matrix over $R$ and a $n \times m$ matrix over $M$ we define 
$A \otimes B$ to be the $m \times m $ matrix with entries in $R \otimes M$ given by 
$(A \otimes B)_{ij} = \displaystyle \sum_{k=1}^{n} A_{ik} \otimes B_{kj} $.
The trace of  $A \otimes B$, written trace$(A \otimes B)$, is given by 
$ \displaystyle \sum_{i=1}^{m} \displaystyle \sum_{k=1}^{n} A_{ik} \otimes B_{ki} 
\in C_{1}(R,M) $. We have that the $1-$chain trace$(A \otimes B)$ is a cycle if and only if 
trace$(AB) = $ trace$(B\phi(A))$, in which case we denote its homology class by $T_{1}(A \otimes B) 
\in HH_{1}(R,M)$.

\subsection{One-parameter Fixed Point Theory}

Let $X$ be a finite connected CW complex and $F: X \times I \to X$ a cellular homotopy.
We consider  $I = [0, 1]$ with the usual CW structure and orientation of cells, 
and $X \times I$ with the product CW structure, where its cells are given the product orientation.

Pick a basepoint $(v,0) \in X \times I$, and a basepath $\tau$ in $X$ from $v$ to $F(v,0)$. 
We identify $\pi_{1}(X \times I, (v,0)) \equiv G$ with $\pi_{1}(X,v) $ via the isomorphism 
induced by projection $p: X \times I \to X$. We write $\phi: G \to G$ for the homomorphism;
$$ \pi_{1}(X \times I, (v,0)) \stackrel{F_{\#}}{\to} \pi_{1}(X, F(v,0)) \stackrel{c_{\tau}}{\to} \pi_{1}(X, v) $$

We choose a lift $\tilde{E}$ in the universal cover, $\tilde{X}$, of $X$ for each cell $E$ and 
we orient $\tilde{E}$ compatibly with $E$. 
Let $\tilde{\tau}$ be the lift of the basepath $\tau$ which starts at in the basepoint $\tilde{v} \in \tilde{X}$ 
and $\tilde{F}: \tilde{X} \times I \to \tilde{X}$ the unique lift of $F$ satisfying $\tilde{F}(\tilde{v},0) = \tilde{\tau}(1)$.

We can regard $C_{\ast}(\tilde{X})$ as a right $\mathbb{Z}G$ chain complex as follows: if $\omega$ is a loop 
at $v$ which lifts to a path $\tilde{\omega}$ starting at $\tilde{v}$ then $\tilde{E}[\omega]^{-1} = h_{[w]}(\tilde{E})$, 
where $h_{[\omega]} $ is the covering transformation sending $\tilde{v}$ to $\tilde{\omega}(1)$.

The homotopy $\tilde{F}$ induces a chain homotopy $\tilde{D_{k}}: C_{k}(\tilde{X}) \to C_{k+1}(\tilde{X})$ 
given by $\tilde{D_{k}}(\tilde{E}) = (-1)^{k+1}F_{k}(\tilde{E} \times I) \in C_{k+1}(\tilde{X})$, for each cell $\tilde{E} \in \tilde{X}$. 
This chain homotopy satisfies; $\tilde{D}(\tilde{E}g) = \tilde{D}(\tilde{E}) \phi(g)$ and the boundary operator 
$\tilde{\partial_{k}}: C_{k}(\tilde{X}) \to C_{k-1}(\tilde{X})$ satisfies; $\tilde{\partial}(\tilde{E}g)= \tilde{\partial}(\tilde{E})g$.

Define endomorphism of, $\oplus_{k} C_{k}(\tilde{X})$, by $\tilde{D_{\ast}} = \oplus_{k} (-1)^{k+1} \tilde{D_{k}}$,
$ \tilde{\partial_{\ast}} = \oplus_{k} \tilde{\partial_{k}}$, $\tilde{F_{0 \ast}} = \oplus_{k} (-1)^{k} \tilde{F_{0 k}} $ 
and $\tilde{F_{1 \ast}} = \oplus_{k} (-1)^{k} \tilde{F_{1 k}} $.
We consider  trace$(\tilde{\partial_{\ast}} \otimes \tilde{D_{\ast}}) \in HH_{1}(\mathbb{Z}G, (\mathbb{Z}G)^{\phi})$. 
This is a Hochshcild 1-chain whose boundary is:
trace$(\tilde{D_{\ast}}\phi(\tilde{\partial_{\ast}}) - \tilde{\partial_{\ast}} \tilde{D_{\ast}}) .  $

We denote by $G_{\phi}(\partial(F))$ the subset of $G_{\phi}$ consisting of semiconjugacy classes associated to fixed 
points of $F_{0}$ or $F_{1}$.

\begin{definition}
The one-parameter trace of homotopy $F$ is;
$$R(F) \equiv T_{1}(\tilde{\partial_{\ast}} \otimes \tilde{D_{\ast}}; G_{\phi}(\partial(F))) \in 
\bigoplus_{C \in G_{\phi} - G_{\phi}(\partial(F))} HH_{1}(\mathbb{Z}G, (\mathbb{Z}G)^{\phi})_{C}  $$
$$ \cong \bigoplus_{C \in G_{\phi} - G_{\phi}(\partial(F))} H_{1}(Z(g_{C})). $$
\end{definition}

{\definition
The C-component of $R(F)$ is denoted by $i(F,C) \in {HH_{1}(\mathbb{Z}G, (\mathbb{Z}G)^{\phi})}_{C}.$
We call it the {\it fixed point index of $F$} corresponding to semiconjugacy class $C \in G_{\phi}$. 
The one-parameter Nielsen number, $N(F)$, of $F$ is the number of nonzero fixed point indices.}
  
\bigskip 
  
The one-parameter Lefschetz class, $L(F)$, of $F$ is defined by;  
$$L(F) = \displaystyle \sum_{C \in G_{\phi} - G_{\phi}(\partial F)} j_{C}(i(F,C)) $$ where 
$j_{C}: H_{1}(Z(g_{C})) \to H_{1}(G)$ is induced by the inclusion $Z(g_{C}) \subset G$.
From \cite{G-N-94} we have the following theorems:

\begin{theorem}[one-parameter Lefschetz fixed point theorem]
If $L(F) \neq 0$ then every map homotopic to $F$ relative to $X \times \{0,1\}$ has a fixed point 
not in the same fixed point class as any fixed point in $X \times \{0,1\}$. In particular, if 
$F_{0}$ and $F_{1}$ are fixed point free, every map homotopic to $F$ relative to $X \times \{0,1\}$ has 
a fixed point.   
\end{theorem}   

\begin{theorem}[one-parameter Nielsen fixed point theorem]
Every map homotopic to $F$ relative to $X \times \{0,1\}$ has at least $N(F)$ fixed point classes 
other than the fixed point classes which meet $X \times \{0,1\}$. In particular, if $F_{0}$ and $F_{1}$ 
are fixed point free maps, then every map homotopic to $F$ relative to $X \times \{0,1\}$ has 
at least $N(F)$ path components.
\end{theorem}

\subsection{Semiconjugacy classes in the torus}

In this subsection we describe some results about the semiconjugacy classes in the torus.

We take $w=[(0,0)] \in T$ and $G = \pi_{1}(T,w) = \{ u, v| uvu^{-1}v^{-1}=1 \}$, where 
$u \equiv a$ and $v \equiv b$. Thus, given a homomorphism $\phi: G \to G$ we have 
$\phi(u) = u^{b_{1}} v^{b_{2}}$ and $\phi(v) = u^{b_{3}} v^{b_{4}}$. Therefore, 
$\phi(u^{m}v^{n}) = u^{mb_{1}+nb_{3}} v^{mb_{2}+nb_{4}}$, for all $m,n \in \mathbb{Z}$.
We denote this homomorphism by the matrix;
$$ [\phi] = \left ( 
\begin{array}{cc}
  b_{1} & b_{3} \\
  b_{2} & b_{4} \\
\end{array} \right )
$$
\begin{prop}
Two elements  $ g_{1} = u^{m_{1}}v^{n_{1}}$ and  $ g_{2} = u^{m_{2}}v^{n_{2}} $ in $G$ belong to the same conjugacy class,
 if and only if there are integers $m,n$ satisfying the following equations:
$$ \left \{
\begin{array}{c}
 m(b_{1}-1)+nb_{3} = m_{2} - m_{1} \\
 mb_{2}+n(b_{4}-1) = n_{2} - n_{1} \\
 \end{array} \right. 
$$ 
\end{prop}
{\it Proof.}
If there is $g = u^{m} v^{n} \in G$ satisfying $g_{1} =g g_{2} \phi(g)^{-1}$ then we obtain the equation of the proposition. 
The other direction is analogous.
\qed

We take the isomorphism $\Theta : G \to \mathbb{Z} \times \mathbb{Z}$ such that $\Theta(u^{m}v^{n}) = (m, n)$.
By above proposition two elements $ g_{1} = u^{m_{1}}v^{n_{1}}$ and  $ g_{2} = u^{m_{2}}v^{n_{2}} $ in $G$
 belong to the same conjugacy class, if and only if 
there is $z \in \mathbb{Z} \times \mathbb{Z}$ satisfying;  
$([\phi]-I)z = \Theta(g_{2}g_{1}^{-1})$, where 
$I$ is the identity matrix. If $det([\phi]-I) \neq 0 $ will have an infinite amount of semiconjugacy classes.

\begin{corollary}
The semicentralizer $Z(g)$ of a element $g \in G$ is isomorphic to the kernel of $[\phi]-I$.
\end{corollary}

\begin{lem}
The 1-chain,  $u^{k} v^{l} \otimes u^{m} v^{n}$, is a cycle if and only if the element $(k,l) \in \mathbb{Z} \times \mathbb{Z}$ 
belongs to the kernel of $[\phi]-I$.
\end{lem}
{\it Proof.}
If $u^{k} v^{l} \otimes u^{m} v^{n} $ is a cycle, then 
$0 = d_{1}( u^{k} v^{l} \otimes u^{m} v^{n} ) = $
$ u^{m} v^{n} \phi(u^{k} v^{l}) - u^{k} v^{l} u^{m} v^{n} = $ \break
$ u^{m} v^{n} u^{kb_{1}+lb_{3}} v^{kb_{2}+lb_{4}} - u^{k} v^{l} u^{m} v^{n} = $ 
$u^{m + kb_{1}+lb_{3}} v^{kb_{2}+lb_{4}+n} - u^{k+m} v^{l+n} $. This implies 
$k(b_{1}-1)+lb_{3} = 0$ and $kb_{2}+l(b_{4}-1) = 0$. 
The other direction is analogous. \qed

\bigskip

\begin{prop}
The 1-chain, $u^{k} \otimes u^{m} v^{n}$, is homologous to the 1-chain, 
$k u \otimes u^{m+k-1} v^{n} $, for all $k,m,n \in \mathbb{Z}$. 
\end{prop}
{\it Proof.}
Note that for $k= 0$ and $1$ the proposition is true. We suppose that  
for some $s > 0 \in \mathbb{Z}$, $u^{s} \otimes u^{m} v^{n}$  
$ \sim s u \otimes u^{m+s-1} v^{n} $. Considering the to 2-chain
 $u^{s} \otimes u \otimes u^{m}v^{n} $ then we have  
$$
\begin{array}{lll}
 d_{2}(u^{s} \otimes u \otimes u^{m}v^{n} ) 
 & = &  u \otimes u^{m+s}v^{n} - 
 u^{s+1} \otimes  u^{m}v^{n} + u^{s} \otimes  u^{1+m}v^{n} \\
 & \sim &  u \otimes u^{m+s}v^{n} - 
 u^{s+1} \otimes  u^{m}v^{n} + su \otimes  u^{1+m+s-1}v^{n} \\
 & = & (s+1)u \otimes  u^{m+(s+1)-1}v^{n} - u^{s+1} \otimes  u^{m}v^{n}.  \\
\end{array}
$$
Therefore $(s+1)u \otimes  u^{m+(s+1)-1}v^{n} \sim u^{s+1} \otimes  u^{m}v^{n}.$  Using induction, we have the result. 
The case in which $k < 0$ is analogous. \qed

\bigskip

\begin{lem}
Each 1-chain, $\displaystyle \sum^{t}_{i=1} a_{i} u^{k_{i}} v^{l_{i}} \otimes u^{m_{i}} v^{n_{i}}  $, 
is homologous to a 1-chain, $\displaystyle \sum^{\bar{t}}_{i=1} \bar{a}_{i} u^{\bar{k_{i}}} v^{\bar{l_{i}}}\otimes 
u^{\bar{m_{i}}} v^{\bar{n_{i}}}  $, where all elements $\bar{l_{i}},$ $i=1,...,\bar{t},$ are positive.
\end{lem}
{\it Proof.}
We denote by $w_{i} = a_{i} u^{k_{i}} v^{l_{i}} \otimes u^{m_{i}} v^{n_{i}}$ and 
$\alpha = \displaystyle \sum^{t}_{i} a_{i} u^{k_{i}} v^{l_{i}} \otimes u^{m_{i}} v^{n_{i}}$.
If there is some $l_{i} \leq 0$ then considering the to 2-chain $\gamma_{i} = a_{i}u^{k_{i}} v^{l_{i}} 
\otimes u^{k_{i}} v^{-l_{i}} \otimes u^{m_{i}-k_{i}} v^{n_{i}-l_{i}} $  
we obtain; $d_{2}(\gamma_{i})= w_{i} - g_{i}+h_{i}$, where
$g_{i} = -a_{i} u^{2k_{i}} \otimes u^{m_{i}-k_{i}} v^{n_{i}-l_{i}} $ and 
$h_{i} = a_{i} u^{k_{i}} v^{-l_{i}} \otimes u^{m_{i}+k_{i}(b_{1}-1)+l_{i}b_{3}} $  $ v^{n_{i}+k_{i}b_{2}+l_{i}(b_{4}-1)} $.
Thus, $w_{i} \sim g_{i} - h_{i}$, and $g_{i}$, $h_{i}$ have the desired form. \qed

\bigskip
In the following proposition we consider $b_{1} = 1$ and $b_{2} = 0$.

\begin{prop}
If the Hochschild 1-chain; $\displaystyle \sum^{t}_{i=1} a_{i} u^{k_{i}} v^{l_{i}} \otimes u^{m_{i}} v^{n_{i}}  $, 
is a 1-cycle then the 1-chain ; $\displaystyle \sum^{t}_{i=1} a_{i} u^{k_{1}+...+k_{t}} v^{l_{1}+...+l_{t}} \otimes u^{m} v^{n}  $, 
is a 1-cycle for all $m,n \in \mathbb{Z}$. 
\end{prop}
{\it Proof.} 
We take a 1-chain,  
$  \displaystyle \sum^{t}_{i} a_{i} u^{k_{i}} v^{l_{i}} \otimes u^{m_{i}} v^{n_{i}}  $, with 
$ d_{1} (\displaystyle \sum^{t}_{i} a_{i} u^{k_{i}} v^{l_{i}} 
 \otimes u^{m_{i}} v^{n_{i}}) =  $ \break $ \displaystyle \sum^{t}_{i} a_{i} u^{m_{i}+k_{i}b_{1}+l_{i}b_{3}} v^{l_{i}b_{4}+k_{i}b_{2}+n_{i}} 
 -a_{i} u^{m_{i}+k_{i}} v^{l_{i}+n_{i}} =0$.
We denote   $e_{i} = $  $u^{m_{i}+k_{i}b_{1}+l_{i}b_{3}} v^{l_{i}b_{4}+k_{i}b_{2}+n_{i}} $ and  
$ f_{i} = u^{m_{i}+k_{i}} v^{l_{i}+n_{i}} $. The last equality implies the following equality on the ring group  $\mathbb{Z}G$:
$$ \displaystyle \sum^{t}_{i} a_{i} e_{i} = \displaystyle \sum^{t}_{i} a_{i} f_{i} .$$
Thus, for each $i$, $1 \leq i \leq t$ there is $j$, $ 1 \leq j \leq t$ such that 
$a_{i}=a_{j}$ and $e_{i} = f_{j}$, that is, we have
$$   (I) \,\, \left \{
\begin{array}{c}
  m_{i} + k_{i}b_{1} +l_{i}b_{3} = k_{j} +m_{j} \\
  l_{i}b_{4} + k_{i}b_{2} + n_{i} = l_{j} + n_{j} \\
 \end{array} \right. 
$$ 

If $i=j $ then the equality above says that the vector $(k_{i},l_{i})$  
satisfies the equation; $([\phi]-I)\left(^{k_{i}}_{l_{i}} \right) = 0$, i.e , 
belongs the to kernel of the $[\phi]-I$. If $i \neq j$ then fixing $j$  there is $q$, $ 1 \leq q \leq t$
such that $a_{j} = a_{q}$ and $e_{j} = f_{q}$. This implies the following equation: 
$$
(II) \,\, \left \{
 \begin{array}{c}
  m_{j} + k_{j}b_{1} +l_{j}b_{3} = k_{q} +m_{q} \\
  l_{j}b_{4} + k_{j}b_{2}+n_{j} = l_{q} + n_{q} \\
 \end{array} \right. 
  $$
Adding the corresponding lines of the $(I)$  and $(II)$ we obtain; 
$$
 \left \{
 \begin{array}{c}
 (k_{i}+k_{j})(b_{1}-1)+(l_{i}+l_{j})b_{3} = k_{q}-k_{i}+m_{q}-m_{i} \\
  (k_{i}+k_{j})b_{2}+ (l_{i}+l_{j})(b_{4}-1) = l_{q}-l_{i}+n_{q}-n_{i}\\
 \end{array} \right. 
  $$ 
If $i=q$ then $(k_{i}+k_{j})(b_{1}-1)+(l_{i}+l_{j})b_{3} = 0$ and 
$(k_{i}+k_{j})b_{2}+ (l_{i}+l_{j})(b_{4}-1) = 0$, which is equivalent to say that the vector, $(k_{i}+k_{j},l_{i}+l_{j})$, 
satisfies the equation; $([\phi]-I)\left(^{k_{i}+k_{j}}_{l_{i}+l_{j}} \right) = 0$. 
Thus, we can take a new $i^{'}$, $1\leq i^{'} \leq t$, and do the same process. 
If $i \neq q$ then we can do the same process above and to obtain a new equation, $(III)$, exactly like in the equation $(II)$, 
and so forth.

Therefore, after making the process for all indices $1 \leq i \leq t $, just add 
all vectors and conclude that the vector; $(\displaystyle \sum^{t}_{j} k_{j}, \displaystyle \sum^{t}_{j} l_{j})$ 
belongs to the kernel of the $([\phi]-I)$. Thus, the 1-chain $\displaystyle \sum^{t}_{i=1} a_{i} u^{k_{1}+...+k_{t}} v^{l_{1}+...+l_{t}} \otimes u^{m} v^{n}  $,
is a cycle, for all $m,n \in \mathbb{Z}$. \qed

\bigskip

Note that if the homomorphism $\phi$ is induced by a homotopy which is induced by a fiber-preserving map as in 
Theorem \ref{daci1theorem}, then the set $\{u \otimes u^{m} v^{n} | m,n \in \mathbb{Z}\}$ 
is one generating set for $ HH_{1}(\mathbb{Z}G,(\mathbb{Z}G)^{\phi})$ $ \cong \bigoplus_{C \in G_{\phi}} H_{1}(Z(g_{C})) $ 
$ \cong \bigoplus_{C \in G_{\phi}} \mathbb{Z} $. 
Since $u \otimes u^{m-1}v^{n} $ $\sim u^{-1} \otimes u^{m+1}v^{n} $, for all $m,n \in \mathbb{Z}$  
we use the generating set $\{ u^{-1} \otimes u^{m} v^{n}| m,n \in \mathbb{Z} \}$.

\section{Computing the number $MF_{S^{1}}[f]$} \label{section-4}

In this section we prove the following theorem:

\begin{theorem}[Main theorem] \label{maintheorem}
If $f: M(A) \to M(A)$ is a fiber-preserving map then the homomorphism $f_{\#}: \pi_{1}(M(A)) \to \pi_{1}(M(A))$ is given by; 
$f_{\#}(a) = a $, $f_{\#}(b) = a^{b_{3}} b^{b_{4}} $ and $f_{\#}(c) = a^{c_{1}} b^{c_{2}} c $ where $a,b,c$ are 
generators of $\pi_{1}(M(A),0)$ previously described. If $M(A)$ is one of the fiber bundle given below;

\noindent In the case $II$ $$ A = \left (
\begin{array}{cc}
  1 & 0 \\
  0 & 1 \\
 \end{array} \right ) \,\, and
\,\, B = \left (
\begin{array}{cc}
  1 & n(b_{4}-1) \\
  0 & b_{4} \\
 \end{array} \right ) \,\,\,\,  $$
$$ A = \left (
\begin{array}{cc}
  1 & 0 \\
  0 & 1 \\
 \end{array} \right ) \,\, and
\,\, B = \left (
\begin{array}{cc}
  1 & b_{3} \\
  0 & -1 \\
 \end{array} \right )
 $$
In the case $III$  
$$ A = \left (
\begin{array}{cc}
  1 & 2k \\
  0 & -1 \\
 \end{array} \right ) \,\, and
\,\, B = \left (
\begin{array}{cc}
  1 & b_{3} \\
  0 & b_{4} \\
 \end{array} \right ) \,\,\,\,  $$
$$ 
A = \left (
\begin{array}{cc}
  1 & a_{3} \\
  0 & -1 \\
 \end{array} \right ) \,\, and
\,\, B = \left (
\begin{array}{cc}
  1 & a_{3} \\
  0 & -1 \\
 \end{array} \right ) 
 $$ 
where $n, k, b_{3}, b_{4}, c_{1}, c_{2}, a_{3}  \in \mathbb{Z}$, then the minimal fixed point set of $f$ 
is composed by $|c_{1}(b_{4}-1) -c_{2}b_{3}|$ disjoint circles. 
This implies; $MF_{S^{1}}[f]  = |c_{1}(b_{4}-1) -c_{2}b_{3}| $.
\end{theorem}

\bigskip
Given a fiber-preserving map $f^{'}: M(A) \to M(A)$ fiberwise homotopic to $f$ then the set $Fix(f^{'})$ is composed by circles. 
The phrase \lq\lq minimal fixed point set of $f$ \rq\rq \, in the theorem above means that we consider the minimum in terms of the 
first homology group, that is, we consider; $min \{rank(H_{1}(Fix(f^{'}))) | f^{'} \sim_{B} f  \}$.

\bigskip

{\it Proof.}
Initially let us consider the case $b_{3} = 0$ and $b_{4}-1 \neq 0$.  
In this situation we must have $a_{3} = 0$ in both of cases.
We take the homotopy $F: T \times I \to T$ defined by; 
$$ F([(x,y)],t)  =  \left \{
\begin{array}{lll}
[( x + 2c_{1}t -\frac{1}{2}, b_{4}y ) ]  & if &  0 \leq t \leq \frac{1}{2} \\

[( x + \frac{2c_{1}-1}{2},  b_{4}y + 2c_{2}t -c_{2} )] & if & \frac{1}{2} \leq t \leq 1 \\
\end{array} \right. 
$$

The homotopy $F$ induces a fiber-preserving map $f: M(A) \to M(A)$ defined by $f(<[(x,y)],t>) = <F([(x,y)],t),t>$.  
Note that the induced homomorphism  by $f$ satisfies; $f_{\#}(a) = a $,
$f_{\#}(b) = b^{b_{4}} $ and $f_{\#}(c) = a^{c_{1}} b^{c_{2}} c $. 
The map $f$ has not fixed points in $t = 0,1$. This implies that $Fix(f) \approx Fix(F)$.

We use the one-parameter trace of $F$, $R(F)$, to compute the minimum number $MF_{S^{1}}[f]$. 

We choose the cellular decomposition for $T$ which consist of two 0-cells; 
$E_{1}^{0} = \{[(0,0)]\}$,  $E_{2}^{0} = \{[(\frac{1}{2},0)]\}$, four 1-cells; 
$E_{1}^{1} = \{[(x,0)]| 0 \leq x \leq \frac{1}{2} \}$, $E_{2}^{1} = \{[(x,0)]| \frac{1}{2} \leq x \leq 1 \}$,
$E_{3}^{1} = \{[(0,y)]| 0 \leq y \leq 1 \}$, $E_{4}^{1} = \{[(\frac{1}{2},y)]| 0 \leq y \leq 1 \}$ and 
two 2-cells; $E_{1}^{2} = \{[(x,y)]| 0 \leq x \leq \frac{1}{2},  0 \leq y \leq 1 \}$, 
$E_{2}^{2} = \{[(x,y)]| \frac{1}{2} \leq x \leq 1,   0 \leq y \leq 1 \}$. For this decomposition the homotopy $F$ is cellular.
 
\begin{figure}[!htp]
\centering \includegraphics[scale=0.25]{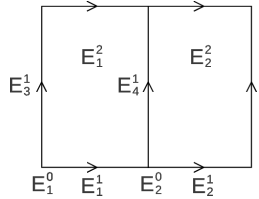}
\caption{Cellular decomposition, case $b_{4}-1 \neq 0$ and $b_{3}= 0$.}
\label{Figura12}
\end{figure}

We orient the cells above as in Figure \ref{Figura12}. By Proposition $4.1$ of \cite{G-N-94} the one-parameter 
trace $R(F)$ is independent of the choice of orientation of cells and the choice of 
lifts to the universal cover. 

For cellular decomposition above we choose in the universal cover $\mathbb{R}^{2}$ 
the lifts which consist of two 0-cells; 
$\tilde{E}_{1}^{0} = (0,0)$,  $\tilde{E}_{2}^{0} = (\frac{1}{2},0)$, four 1-cells; 
$\tilde{E}_{1}^{1} = \{(x,0)| 0 \leq x \leq \frac{1}{2} \}$, $\tilde{E}_{2}^{1} = \{(x,0)| \frac{1}{2} \leq x \leq 1 \}$,
$\tilde{E}_{3}^{1} = \{(0,y)| 0 \leq y \leq 1 \}$, $\tilde{E}_{4}^{1} = \{(\frac{1}{2},y)| 0 \leq y \leq 1 \}$ and  
two 2-cells; $\tilde{E}_{1}^{2} = \{(x,y)| 0 \leq x \leq \frac{1}{2},  0 \leq y \leq 1 \}$, 
$\tilde{E}_{2}^{2} = \{(x,y)| \frac{1}{2} \leq x \leq 1,   0 \leq y \leq 1 \}$.

We consider $w = [(0,0)]$ the basepoint and $\tau$ basepath, the linear path between $w$ and 
$F(w,0)$. We take the lifts $\tilde{w}= (0,0)$ and $\tilde{\tau}$ the linear path between $\tilde{w}$ and $(-\frac{1}{2},0)$.  The unique lift
$\tilde{F}: \mathbb{R}^{2} \times I \to \mathbb{R}^{2} $ of $F$ mapping $(\tilde{w},0)$ to $\tilde{\tau}(1)$  
 is given by;
$$ \tilde{F}(x,y,t) = \left \{
\begin{array}{lll}
( x + 2c_{1}t -\frac{1}{2}, b_{4}y )   & if & 0 \leq t \leq \frac{1}{2} \\

( x + \frac{2c_{1}-1}{2},  b_{4}y + 2c_{2}t -c_{2} ) & if & \frac{1}{2} \leq t \leq 1 \\
\end{array} \right. 
$$

If $G = \pi_{1}(T,[(0,0)]) = \{ u,v| uvu^{-1}v^{-1} = 1\}$ then matrices of operators $\tilde{\partial}_{1}$, $\tilde{\partial}_{2}$, 
$\tilde{D}_{0}$ and $\tilde{D}_{1}$ are given by;

$$ [\tilde{\partial}_{1}] = \left (
\begin{array}{cccc}
  -1 & u^{-1} & v^{-1}-1 & 0 \\
   1 & -1 & 0 & v^{-1}-1  \\
\end{array} \right ) 
$$

$$ [\tilde{\partial}_{2}] = \left (
\begin{array}{cc}
  v^{-1}-1 & 0 \\
  0 & v^{-1}-1 \\
  1 & -u^{-1} \\
  -1 & 1 \\
\end{array} \right ) 
$$

$$ [\tilde{D}_{0}] = \left (
\begin{array}{cc}
  -\tilde{X}(c_{1}) & -\tilde{X}(c_{1}) \\
 -\tilde{Y}(c_{1}) & -\tilde{X}(c_{1}) \\
  0 & -u^{-c_{1}} \tilde{W}(c_{2}) \\
 -u^{1-c_{1}} \tilde{W}(c_{2}) & 0 \\
\end{array} \right ) 
$$

$$ [\tilde{D}_{1}] = \left (
\begin{array}{cccc}
 0  & u^{1-c_{1}}\tilde{W}(c_{2}) &  \tilde{X}(c_{1}) \tilde{W}(b_{4}) & \tilde{X}(c_{1}) \tilde{W}(b_{4}) \\
  u^{1-c_{1}}\tilde{W}(c_{2}) & 0 &  \tilde{Y}(c_{1}) \tilde{W}(b_{4}) & \tilde{X}(c_{1}) \tilde{W}(b_{4}) \\
 \end{array} \right ) 
$$
where
{ \small
$$ \tilde{X}(m) = \left \{
\begin{array}{cc}
  \displaystyle\sum^{m}_{j=1} u^{1-j}   & if \,\,\,\,  m > 0   \\
  0 & if \,\,\,\, m=0 \\
 \displaystyle\sum^{-m}_{j=1} -u^{j}   & if \,\,\,\, m < 0   \\
 \end{array} \right. ,
 \hspace{1.0cm}  \tilde{Y}(m) = \left \{
\begin{array}{cc}
  \displaystyle\sum^{m}_{j=1} u^{2-j}   & if \,\,\,\,  m > 0   \\
  0 & if \,\,\,\, m=0 \\
 \displaystyle\sum^{-m}_{j=1} -u^{j+2}   & if \,\,\,\, m < 0   \\
 \end{array} \right. $$ }

{  \small
$$ and \hspace{2.0cm}  \tilde{W}(m) = \left \{
\begin{array}{cc}
   \displaystyle\sum^{m}_{j=1} v^{1-j}  & if \,\,\,\,  m > 0   \\
  0  & if \,\,\,\, m=0 \\
 \displaystyle\sum^{-m}_{j=1} -v^{j}   & if \,\,\,\, m < 0   \\
 \end{array} . \right. $$ }

Thus we have; 
$$
R(F) = T_{1}(\tilde{\partial_{\ast}} \otimes \tilde{D_{\ast}}) = u^{-1}\otimes \tilde{Y}(c_{1}) -2 \otimes \tilde{X}(c_{1})  
 +  1 \otimes \tilde{Y}(c_{1}) $$ $$+  1 \otimes \tilde{X}(c_{1}) \tilde{W}(b_{4}) -u^{-1} \otimes \tilde{Y}(c_{1}) \tilde{W}(b_{4}).
$$

Two elements $g_{1},g_{2} \in G$ belong to the same conjugacy class if and only if there is an 
element $g \in G$ satisfying to equation: $g_{1} = g g_{2} \phi(g^{-1})$.
In this case two elements $u^{-1} \otimes u^{m}v^{s}$ and $u^{-1} \otimes u^{n}v^{t}$, $m, n, s, t \in \mathbb{Z},$ 
belong the same semiconjugate class if and only if there is $k \in \mathbb{Z}$ satisfying;
$$ \left \{
\begin{array}{l}
 m = n \\
 s = k(b_{4}-1)+t \\
\end{array} \right.
$$

\begin{figure}[!htp]
\centering
 \includegraphics[scale=0.3]{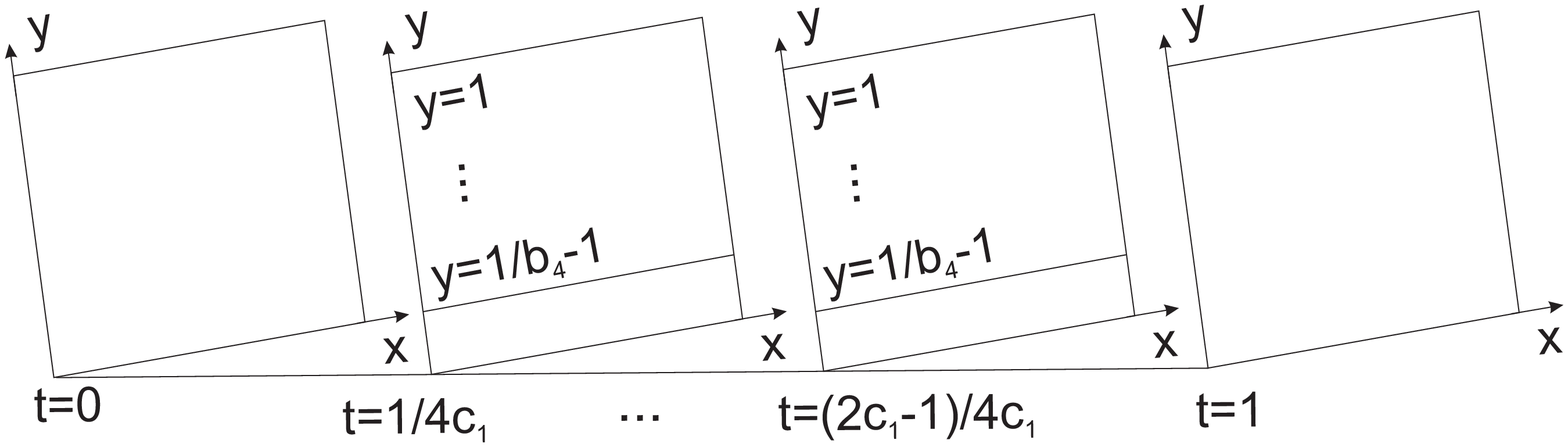}
 \caption{The set Fix(F) in the case $c_{1}$ and $ b_{4}-1$ positives.}
 \label{Figura1}
\end{figure}

If $c_{1}(b_{4}-1) \neq 0$ then we have $N(F) = |c_{1}(b_{4}-1)|$. Since $Fix(F)$ consist of 
$|c_{1}(b_{4}-1)|$ circles, $MF[F]$ is composed of $|c_{1}(b_{4}-1)|$ disjoint circles, see Figure \ref{Figura1}.
Therefore the minimal fixed point set of $f$ consist of $|c_{1}(b_{4}-1)|$ disjoint circles.
If $c_{1}(b_{4}-1) = 0$ then $MF_{S^{1}}[f] = 0$. 
Thus the number $MF_{S^{1}}[f] = |c_{1}(b_{4}-1)|$.

\vspace{0.3cm}

Now in the case II with $b_{3} = n(b_{4}-1)$, $n \in \mathbb{Z}$ we take the homotopy 
$F: T \to T$ given by $F([(x,y)],t) = [(x+b_{3}y +c_{1}t- \frac{1}{2}, b_{4}y +c_{2}t)]$ and 
the fiber-preserving map $f: M(A) \to M(A)$ induced by $F$. We consider the isomorphism of fiber bundle 
$P: M(A) \to M(A^{1}) $,  $A^{1} = P \circ A \circ P^{-1}$, induced by the isomorphism on torus which 
also is denoted by $P: T \to T$ given by the following matrix:
$$ P = \left (
\begin{array}{cc}
  1 & -n \\
  0 & 1
\end{array} \right )
$$

By propoposition \ref{propconjugation} the fiber-preserving map $g: M(A^{1}) \to M(A^{1})$ induced by homotopy 
$G = P \circ F \circ (P^{-1} \times I)$ has $MF_{S^{1}}[g] = MF_{S^{1}}[f]$. Here the homotopy $G$  
is given by; $G([(x,y)],t) =[( x+ (c_{1} -nc_{2})t -\frac{1}{2}, b_{4}y + c_{2}t )] $. 
$$
\xymatrix{ M(A)  \ar[r]^{f} \ar[d]_{P} &   M(A) \ar[d]^{P}  \\
M(A^{1}) \ar[r]_{g}  &   M(A^{1}) \\ }
$$

Note that the homotopy $G$ is homotopic, relative to $T \times \{0,1\}$, to the 
homotopy $G^{'}$ given by;
$$ G^{'}([(x,y)],t) = \left \{
\begin{array}{ll}
[( x +  2(c_{1}-nc_{2})t - \frac{1}{2} , b_{4}y ) ]  & ,  0 \leq t \leq \frac{1}{2} \\

[( x + (c_{1}-nc_{2}) - \frac{1}{2} , b_{4}y  + 2c_{2}t -c_{2} )] & ,  \frac{1}{2} \leq t \leq 1 \\
\end{array} \right. 
$$ 
In fact, using the notation $G([(x,y)],t) = [(\alpha(x,y,t), \beta(x,y,t))]$, where   
$\alpha(x,y,t) = x + (c_{1}-nc_{2})t - \frac{1}{2}$ and $ \beta(x,y,t)  = b_{4}y +c_{2}t $, 
then $H: T \times I \times I \to T $  
defined by;
$$ H([(x,y)],t,s) = \left \{
\begin{array}{lll}
 [(\alpha(x,y,t), \beta(x,y,t))] & if & 0 \leq t \leq s \\
 
 [(\alpha(x,y,2t-s), \beta(x,y,s))] & if & s \leq t \leq \frac{(1+s)}{2} \\
 
 [(\alpha(x,y,1), \beta(x,y,2t-1))] & if & \frac{(1+s)}{2} \leq t \leq 1 \\
\end{array} \right.
$$
is a homotopy, relative to $T \times \{0,1\}$, between $G$ and $G^{'}$. 
Thus, we have $R(G) = R(G^{'})$. Therefore, we can use the previously case and proposition \ref{propconjugation} to show
that the minimal fixed point set of $f$ over $S^{1}$ is composed by  $|c_{1}(b_{4}-1) -c_{2}b_{3}|$ 
disjoint circles.

\bigskip

In the case $III$ we have $ a_{3}(b_{4}-1) = -2b_{3}$. Therefore if $a_{3}$ is even then  
$b_{3} = \frac{-a_{3}}{2}(b_{4} - 1) $. Thus we can use a similar argument 
as in the case above and show that the minimal fixed point set of a fiber-preserving map $f: M(A) \to M(A)$ in this situation 
is composed by $|c_{1}(b_{4}-1) -c_{2}b_{3}|$ disjoint circles.
Note that if $a_{3}$ is even, then a fiber-preserving map $f: M(A) \to M(A)$ in a fiber bundle $M(A)$ with 
$$ A = \left (
\begin{array}{cc}
  1 & a_{3} \\
  0 & -1
\end{array} \right )
$$
has the minimal fixed point set over $S^{1}$ composed by $|c_{1}(b_{4}-1) -c_{2}b_{3}|$ disjoint circles.

\bigskip

Now, let us consider the cases $II$ and $III$ in the following situation; $b_{4} = -1$ and $b_{3} = 2k +1$, $k \in \mathbb{Z}$. 
Note that the case $b_{3}$ even has already been solved. First we take $b_{3}= 1$. 

Consider the fiber-preserving map $f: MA \to MA$ induced by homotopy $F: T \times I \to T$ given by;
$$ F([(x,y)],t)  =  \left \{
\begin{array}{ll}
[( x + y + 2c_{1}t +\frac{1}{2}, -y + \frac{1}{2}) ]  & ,  0 \leq t \leq \frac{1}{2} \\

[( x + y + \frac{2c_{1}+1}{2},  -y + 2c_{2}t -c_{2} + \frac{1}{2})] & ,  \frac{1}{2} \leq t \leq 1 \\
\end{array} \right. 
$$

Note that $f$ has no fixed point in $t=0,1$. For compute the one-parameter trace $R(F)$ we consider the 
cellular decomposition of the torus which consist of four 0-cells; 
$E_{1}^{0} = \{[(0,0)]\}$,  $E_{2}^{0} = \{[(\frac{1}{2},0)]\}$,
$E_{3}^{0} = \{[(0,\frac{1}{2})]\}$, $E_{4}^{0} = \{[(\frac{1}{2},\frac{1}{2})]\}$,  twelve 1-cells; 
$E_{1}^{1} = \{[(x,0)]| 0 \leq x \leq \frac{1}{2} \}$, $E_{2}^{1} = \{[(x,0)]| \frac{1}{2} \leq x \leq 1 \}$,
$E_{3}^{1} = \{[(0,y)]| 0 \leq y \leq \frac{1}{2} \}$, $E_{4}^{1} = \{[(y,-y+ \frac{1}{2})]| 0 \leq y \leq \frac{1}{2} \}$,  
$E_{5}^{1} = \{[(\frac{1}{2},y)]| 0 \leq y \leq \frac{1}{2} \}$,
$E_{6}^{1} = \{[(y,-y+ 1)]| \frac{1}{2} \leq y \leq 1 \}$,
$E_{7}^{1} = \{[(x,\frac{1}{2})]| 0 \leq x \leq \frac{1}{2} \}$,
$E_{8}^{1} = \{[(x,\frac{1}{2})]| \frac{1}{2} \leq x \leq 1 \}$,
$E_{9}^{1} = \{[(0,y)]| \frac{1}{2} \leq y \leq 1 \}$,
$E_{10}^{1} = \{[(y,-y+ 1)]| 0 \leq y \leq \frac{1}{2} \}$,
$E_{11}^{1} = \{[(\frac{1}{2},y)]| \frac{1}{2} \leq y \leq 1 \}$, 
$E_{12}^{1} = \{[(y+\frac{1}{2},-y+1)]| 0 \leq y \leq \frac{1}{2} \}$, and
eight 2-cells; $E_{1}^{2} = \{[(x,y)]| 0 \leq x \leq \frac{1}{2},  0 \leq y \leq -x+\frac{1}{2} \}$, 
$E_{2}^{2} = \{[(x,y)]| 0 \leq x \leq \frac{1}{2},  -x+\frac{1}{2} \leq y \leq \frac{1}{2} \}$, 
$E_{3}^{2} = \{[(x,y)]| \frac{1}{2} \leq x \leq 1,  0 \leq y \leq -x+1 \}$,
$E_{4}^{2} = \{[(x,y)]| \frac{1}{2} \leq x \leq 1,  -x+1 \leq y \leq \frac{1}{2} \}$,
$E_{5}^{2} = \{[(x,y)]| 0 \leq x \leq \frac{1}{2},  \frac{1}{2} \leq y \leq -x+1 \}$,
$E_{6}^{2} = \{[(x,y)]| 0 \leq x \leq \frac{1}{2},  -x+1 \leq y \leq 1 \}$,
$E_{7}^{2} = \{[(x,y)]| \frac{1}{2} \leq x \leq 1,  \frac{1}{2} \leq y \leq -x+\frac{3}{2} \}$,
$E_{8}^{2} = \{[(x,y)]| \frac{1}{2} \leq x \leq 1,  -x+\frac{3}{2} \leq y \leq 1 \}$.

These cells are oriented as in the figure below.
For this cellular decomposition the homotopy $F$ is cellular.

\begin{figure}[htp]
\centering
\includegraphics[scale=0.35]{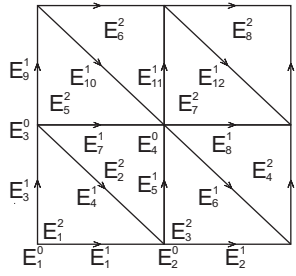}
\caption{Cellular decomposition, case $b_{4}=-1$}
\end{figure}

For the cellular decomposition above we choose in the universal cover $\mathbb{R}^{2}$ 
the lifts which consist of four 0-cells; 
$\tilde{E}_{1}^{0} = \{(0,0)\}$,  $\tilde{E}_{2}^{0} = \{(\frac{1}{2},0)\}$,
$\tilde{E}_{3}^{0} = \{(0,\frac{1}{2})\}$, $\tilde{E}_{4}^{0} = \{(\frac{1}{2},\frac{1}{2})\}$, twelve 1-cells; 
$\tilde{E}_{1}^{1} = \{(x,0)| 0 \leq x \leq \frac{1}{2} \}$, $\tilde{E}_{2}^{1} = \{(x,0)| \frac{1}{2} \leq x \leq 1 \}$,
$\tilde{E}_{3}^{1} = \{(0,y)| 0 \leq y \leq \frac{1}{2} \}$, $\tilde{E}_{4}^{1} = \{(y,-y+ \frac{1}{2})| 0 \leq y \leq \frac{1}{2} \}$,  
$\tilde{E}_{5}^{1} = \{(\frac{1}{2},y)| 0 \leq y \leq \frac{1}{2} \}$,
$\tilde{E}_{6}^{1} = \{(y,-y+ 1)| \frac{1}{2} \leq y \leq 1 \}$,
$\tilde{E}_{7}^{1} = \{(x,\frac{1}{2})| 0 \leq x \leq \frac{1}{2} \}$,
$\tilde{E}_{8}^{1} = \{(x,\frac{1}{2})| \frac{1}{2} \leq x \leq 1 \}$,
$\tilde{E}_{9}^{1} = \{(0,y)| \frac{1}{2} \leq y \leq 1 \}$,
$\tilde{E}_{10}^{1} = \{(y,-y+ 1)| 0 \leq y \leq \frac{1}{2} \}$,
$\tilde{E}_{11}^{1} = \{(\frac{1}{2},y)| \frac{1}{2} \leq y \leq 1 \}$,
$\tilde{E}_{12}^{1} = \{(y+\frac{1}{2},-y+1)| 0 \leq y \leq \frac{1}{2} \}$, and
eight 2-cells; $\tilde{E}_{1}^{2} = \{(x,y)| 0 \leq x \leq \frac{1}{2},  0 \leq y \leq -x+\frac{1}{2} \}$, 
$\tilde{E}_{2}^{2} = \{(x,y)| 0 \leq x \leq \frac{1}{2},  -x+\frac{1}{2} \leq y \leq \frac{1}{2} \}$, 
$\tilde{E}_{3}^{2} = \{(x,y)| \frac{1}{2} \leq x \leq 1,  0 \leq y \leq -x+1 \}$,
$\tilde{E}_{4}^{2} = \{(x,y)| \frac{1}{2} \leq x \leq 1,  -x+1 \leq y \leq \frac{1}{2} \}$,
$\tilde{E}_{5}^{2} = \{(x,y)| 0 \leq x \leq \frac{1}{2},  \frac{1}{2} \leq y \leq -x+1 \}$,
$\tilde{E}_{6}^{2} = \{(x,y)| 0 \leq x \leq \frac{1}{2},  -x+1 \leq y \leq 1 \}$,
$\tilde{E}_{7}^{2} = \{(x,y)| \frac{1}{2} \leq x \leq 1,  \frac{1}{2} \leq y \leq -x+\frac{3}{2} \}$,
$\tilde{E}_{8}^{2} = \{(x,y)| \frac{1}{2} \leq x \leq 1,  -x+\frac{3}{2} \leq y \leq 1 \}$.

\bigskip

We take $w = [(0,0)]$ the basepoint and $\tau$ basepath, the linear path between $w$ and 
$F(w,0)$. We take the lifts $\tilde{w}= (0,0)$ and $\tilde{\tau}$ the linear path between $\tilde{w}$ and $(0,\frac{1}{2})$.  The unique lift
$\tilde{F}: \mathbb{R}^{2} \times I \to \mathbb{R}^{2} $ of $F$ mapping $(\tilde{w},0)$ to $\tilde{\tau}(1)$  
is given by;
$$ \tilde{F}((x,y),t)  =  \left \{
\begin{array}{ll}
( x + y + 2c_{1}t +\frac{1}{2}, -y + \frac{1}{2})   & ,  0 \leq t \leq \frac{1}{2} \\
( x + y + \frac{2c_{1}+1}{2},  -y + 2c_{2}t -c_{2} + \frac{1}{2}) & ,  \frac{1}{2} \leq t \leq 1 \\
\end{array} \right. 
$$

If $G = \pi_{1}(T,[(0,0)]) = \{ u,v| uvu^{-1}v^{-1} = 1\}$ then 
$$ [\tilde{\partial}_{1}] = \left ( 
\begin{array}{cccccccccccc}
 -1 & u^{-1} & -1 & 0 & 0 & u^{-1} & 0 & 0 & v^{-1} & -v^{-1} & 0 & 0 \\
  1 & -1 & 0 & 1 & -1 & 0 & 0 & 0 & 0 & 0 & v^{-1} & -v^{-1}  \\
 0 & 0 & 1 & -1 & 0 & 0 & -1 & u^{-1} & -1  & 0 & 0 & u^{-1} \\
 0 & 0 & 0 & 0 & 1 & -1 & 1 & -1 & 0 & 1 & -1 & 0 \\
\end{array} \right ),
$$

$${ \small [\tilde{\partial}_{2}] = \left ( 
\begin{array}{cccccccc}
 1 & 0 & 0 & 0 & 0 & -v^{-1} & 0 & 0 \\
 0 & 0 & 1 & 0 & 0 & 0 & 0 & -v^{-1} \\
 -1 & 0 & 0 & u^{-1} & 0 & 0 & 0 & 0  \\
 -1 & 1 & 0 & 0 & 0 & 0 & 0 & 0  \\ 
 0 & 1 & -1 & 0 & 0 & 0 & 0 & 0 \\
 0 & 0 & -1 & 1 & 0 & 0 & 0 & 0 \\
 0 & -1 & 0 & 0 & 1 & 0 & 0 & 0 \\
 0 & 0 & 0 & -1 & 0 & 0 & 1 & 0 \\
 0 & 0 & 0 & 0 & -1 & 0 & 0 & u^{-1} \\
 0 & 0 & 0 & 0 & -1 & 1 & 0 & 0 \\
 0 & 0 & 0 & 0 & 0 & 1 & -1 & 0 \\
 0 & 0 & 0 & 0 & 0 & 0 & -1 & 1 \\ 
\end{array} \right ) },
$$

$$ {\small [\tilde{D}_{0}] = \left ( 
\begin{array}{cccc}
 0 & 0 & -u^{-1} \tilde{X}(c_{1}) & -u^{-2} \tilde{X}(c_{1}) \\
 0 & 0 &  -u^{-1} \tilde{X}(c_{1}) & -u^{-1} \tilde{X}(c_{1}) \\
 0 & -v^{-1} \tilde{W}(c_{2}) & -u^{-c_{1}-1} \tilde{W}(c_{2}) & 0 \\
 0 & 0 & 0 & 0 \\
 -u^{-c_{1}}v^{-1} \tilde{W}(c_{2}) & 0 & 0 & -u^{-c_{1}-1} \tilde{W}(c_{2}) \\
 0 & 0 & 0 & 0 \\
 -u^{-c_{1}} \tilde{X}(c_{1}) & -u^{-1} \tilde{X}(c_{1}) & 0 & 0 \\
 -\tilde{X}(c_{1}) & -u^{-1} \tilde{X}(c_{1}) & 0 & 0 \\
 0 & -u^{-c_{1}-1} \tilde{W}(c_{2}) & -u^{-c_{1}-1} \tilde{W}(c_{2}) & 0 \\
 0 & 0 & 0 & 0 \\
 -u^{-c_{1}} \tilde{W}(c_{2}) & 0 & 0 & -u^{-c_{1}} \tilde{W}(c_{2}) \\
 0 & 0 & 0 & 0 \\ 
\end{array} \right ) }
$$
and with the following data; 

$$\tilde{D}_{1}(\tilde{E}^{1}_{1}) = \tilde{E}^{2}_{3}u^{-c_{1}}v^{-1} \tilde{W}(c_{2})  
 + \tilde{E}^{2}_{4}u^{-c_{1}}v^{-1} \tilde{W}(c_{2}) + \tilde{E}^{2}_{7}u^{-c_{1}} \tilde{W}(c_{2}) + 
 \tilde{E}^{2}_{8}u^{-c_{1}} \tilde{W}(c_{2}), $$
 $$\tilde{D}_{1}(\tilde{E}^{1}_{2}) = \tilde{E}^{2}_{1}u^{-c_{1}-1}v^{-1} \tilde{W}(c_{2})  
 + \tilde{E}^{2}_{2}u^{-c_{1}-1}v^{-1} \tilde{W}(c_{2}) + \tilde{E}^{2}_{5}u^{-c_{1}-1} \tilde{W}(c_{2}) + 
 \tilde{E}^{2}_{6}u^{-c_{1}-1} \tilde{W}(c_{2}), $$
$$
\begin{array}{lll}
\tilde{D}_{1}(\tilde{E}^{1}_{3}) & = &
\tilde{E}^{2}_{1}u^{-c_{1}} \tilde{X}(c_{1}) + 
 \tilde{E}^{2}_{2}u^{-c_{1}} \tilde{X}(c_{1})  
 + \tilde{E}^{2}_{3}(u^{-c_{1}}\tilde{X}(c_{1})+u^{-c_{1}}v^{-1} \tilde{W}(c_{2}) ) \\
 & + &   \tilde{E}^{2}_{4}( \tilde{X}(c_{1})+u^{-c_{1}} \tilde{W}(c_{2}) ) + 
 \tilde{E}^{2}_{7}u^{-c_{1}} \tilde{W}(c_{2}) + 
 \tilde{E}^{2}_{8}u^{-c_{1}} \tilde{W}(c_{2}), \\ 
 \end{array}
 $$
$$\tilde{D}_{1}(\tilde{E}^{1}_{4}) = \tilde{E}^{2}_{1}u^{-c_{1}} \tilde{X}(c_{1})  
 + \tilde{E}^{2}_{2}u^{-c_{1}} \tilde{X}(c_{1}) + \tilde{E}^{2}_{3}u^{-c_{1}} \tilde{X}(c_{1}) + 
 \tilde{E}^{2}_{4}u^{-c_{1}} \tilde{X}(c_{1}), $$
 $$
\begin{array}{lll}
\tilde{D}_{1}(\tilde{E}^{1}_{5}) & = &
\tilde{E}^{2}_{1}(u^{-2} \tilde{X}(c_{1})+ u^{-c_{1}-1}v^{-1} \tilde{W}(c_{2}) ) + 
 \tilde{E}^{2}_{2}(u^{-1} \tilde{X}(c_{1})+u^{-c_{1}-1} \tilde{W}(c_{2}) ) \\  
 & + & \tilde{E}^{2}_{3}u^{-1}\tilde{X}(c_{1}) 
  +    \tilde{E}^{2}_{4}u^{-1}\tilde{X}(c_{1}) + 
 \tilde{E}^{2}_{5}u^{-c_{1}-1} \tilde{W}(c_{2}) + 
 \tilde{E}^{2}_{6}u^{-c_{1}-1} \tilde{W}(c_{2}), \\ 
 \end{array}
 $$
 $$\tilde{D}_{1}(\tilde{E}^{1}_{6}) = \tilde{E}^{2}_{1}u^{-2} \tilde{X}(c_{1})  
 + \tilde{E}^{2}_{2}u^{-2} \tilde{X}(c_{1}) + \tilde{E}^{2}_{3}u^{-1} \tilde{X}(c_{1}) + 
 \tilde{E}^{2}_{4}u^{-1} \tilde{X}(c_{1}), $$ 
$$
\begin{array}{lll}
\tilde{D}_{1}(\tilde{E}^{1}_{7}) & = & \tilde{E}^{2}_{1}u^{-c_{1}-1}v^{-1} \tilde{W}(c_{2})  
 + \tilde{E}^{2}_{2}u^{-c_{1}-1}v^{-1} \tilde{W}(c_{2}) + 
 \tilde{E}^{2}_{5}u^{-c_{1}-1}v^{-1} \tilde{W}(c_{2}) \\
  & + &  \tilde{E}^{2}_{6}u^{-c_{1}-1}v^{-1} \tilde{W}(c_{2}), \\
 \end{array} 
 $$
$$\tilde{D}_{1}(\tilde{E}^{1}_{8}) = \tilde{E}^{2}_{1}u^{-c_{1}-1} \tilde{W}(c_{2})  
 + \tilde{E}^{2}_{2}u^{-c_{1}-1} \tilde{W}(c_{2}) + \tilde{E}^{2}_{5}u^{-c_{1}-1} \tilde{W}(c_{2}) + 
 \tilde{E}^{2}_{6}u^{-c_{1}-1} \tilde{W}(c_{2}), $$
$$
\begin{array}{lll}
\tilde{D}_{1}(\tilde{E}^{1}_{9}) & = &
\tilde{E}^{2}_{1}u^{-c_{1}-1} \tilde{W}(c_{2}) + 
 \tilde{E}^{2}_{2}u^{-c_{1}-1} \tilde{W}(c_{2})  
 + \tilde{E}^{2}_{5}(u^{-2}v\tilde{X}(c_{1})+u^{-c_{1}-1}\tilde{W}(c_{2}) ) \\
 & + &   \tilde{E}^{2}_{6}(u^{-1}v \tilde{X}(c_{1})+u^{-c_{1}-1}v \tilde{W}(c_{2}) ) + 
 \tilde{E}^{2}_{7}u^{-1}v \tilde{X}(c_{1}) + 
 \tilde{E}^{2}_{8}u^{-1}v \tilde{X}(c_{1}), \\ 
 \end{array}
 $$
 $$\tilde{D}_{1}(\tilde{E}^{1}_{10}) = \tilde{E}^{2}_{5}u^{-2}v \tilde{X}(c_{1})  
 + \tilde{E}^{2}_{6}u^{-2}v \tilde{X}(c_{1}) + \tilde{E}^{2}_{7}u^{-1}v \tilde{X}(c_{1}) + 
 \tilde{E}^{2}_{8}u^{-1}v \tilde{X}(c_{1}), $$ 
$$
\begin{array}{lll}
\tilde{D}_{1}(\tilde{E}^{1}_{11}) & = &
\tilde{E}^{2}_{3}u^{-c_{1}-1} \tilde{W}(c_{2}) + 
 \tilde{E}^{2}_{4}u^{-c_{1}-1} \tilde{W}(c_{2})  
 + \tilde{E}^{2}_{5} u^{-2}v\tilde{X}(c_{1}) + 
  \tilde{E}^{2}_{6}u^{-2}v \tilde{X}(c_{1}) \\
 & + &    \tilde{E}^{2}_{7}(u^{-2}v \tilde{X}(c_{1}) + u^{-c_{1}-1} \tilde{W}(c_{2} )) + 
 \tilde{E}^{2}_{8}(u^{-1}v \tilde{X}(c_{1}) + u^{-c_{1}-1}v \tilde{W}(c_{2} )), \\ 
 \end{array}
 $$
 $$\tilde{D}_{1}(\tilde{E}^{1}_{12}) = \tilde{E}^{2}_{5}u^{-2}v \tilde{X}(c_{1})  
 + \tilde{E}^{2}_{6}u^{-2}v \tilde{X}(c_{1}) + \tilde{E}^{2}_{7}u^{-2}v \tilde{X}(c_{1}) + 
 \tilde{E}^{2}_{8}u^{-2}v \tilde{X}(c_{1}). $$ 
we can construct the matrix $[\tilde{D}_{1}]_{8 \times 12} $ of the operator $\tilde{D}_{1}$.
Therefore,
$$
\begin{array}{lll}
R(F) & = & -1 \otimes u^{-c_{1}-1} \tilde{W}(c_{2}) -1 \otimes u^{c_{1}} \tilde{W}(c_{2})
-1 \otimes u^{c_{1}} \tilde{X}(c_{1})+ u^{-1} \otimes \tilde{X}(c_{1})\\ 
 & + & u^{-1} \otimes u^{-c_{1}}\tilde{W}(c_{2})
+ 1 \otimes( u^{-1}\tilde{X}(c_{1}) + u^{-c_{1}-1}\tilde{W}(c_{2})) 
 -  1 \otimes u^{-1}\tilde{X}(c_{1}) \\
 & - & 1 \otimes( u^{-2}v\tilde{X}(c_{1}) + u^{-c_{1}-1}\tilde{W}(c_{2})) 
+ u^{-1} \otimes u^{-1}v\tilde{X}(c_{1})  + 1 \otimes u^{-2}v\tilde{X}(c_{1}) \\ 
 & - & 1 \otimes( u^{-2}v\tilde{X}(c_{1}) + u^{-c_{1}-1}\tilde{W}(c_{2})).\\
\end{array}
$$

Similar to the case $b_{3} = 0$ and $b_{4}-1 \neq 0$ we obtain; $$ N(F) = |2c_{1}+c_{2}| = |c_{1}(b_{4}-1)-c_{2}b_{3}|.$$
Since $Fix(F)$ is composed by $|c_{1}(b_{4}-1)-c_{2}b_{3}|$ disjoint circles, then the minimal fixed point set of 
$f: M(A) \to M(A)$ induced by $F: T \times I \to T$ is composed by $|c_{1}(b_{4}-1)-c_{2}b_{3}|$ disjoint circles. 
Therefore, $MF_{S^{1}}[f] = |c_{1}(b_{4}-1)-c_{2}b_{3}|.$

\bigskip

The case $b_{3} = 2k +1$ with $k \neq 0$, we take the fiber-preserving map $f: M(A) \to M(A)$ induced by  
$F: T \times I \to T$ given by $ F([(x,y)],t)  =  [( x + b_{3}y + c_{1}t +\frac{-k+1}{2}, -y + c_{2}t+ \frac{1}{2}) ]  $.
Conjugating the homotopy $F$ by the isomorphism $P: T \to T$ given by 
$$ [P] = \left (
\begin{array}{cc}
 1 & k \\
 0 & 1
\end{array} \right )
$$
we obtain the homotopy $G = P \circ F \circ (P^{-1} \times I)$.
The fiber-preserving map $g: M(A^{1}) \to M(A^{1})$, $A^{1} = P \circ A \circ P^{-1}$, given by $g(<[(x,y)],t>) = <G([(x,y)],t),t>$, 
has $MF_{S^{1}}[g] = MF_{S^{1}}[f]$. By the case above and proposition \ref{propconjugation} we can conclude that the minimal 
fixed point set of $f$ is composed by $|c_{1}(b_{4}-1)-c_{2}b_{3}|$ disjoint circles. This implies 
$MF_{S^{1}}[f] = |c_{1}(b_{4}-1)-c_{2}b_{3}|$.
\qed

 \bigskip

\begin{remark} 
Note that by Theorem \ref{main-theorem-daci1} the number $|c_{1}(b_{4}-1) -c_{2}b_{3}|$ appeared in \cite{G-P-V-04} 
only to decide when a fiber-preserving map, in the cases II and III, can be deformed by a fiberwise homotopy to a 
fixed point free map. In Theorem \ref{maintheorem} we have shown that the number $|c_{1}(b_{4}-1) -c_{2}b_{3}|$ is exactly the number 
of circles of minimal fixed point set of a fiber-preserving map in the cases II and III. 
\end{remark}

{\bf Acknowledgments.} 
I would like to thank Prof. Daniel Vendr\'uscolo and Prof. Jo\~ao P. Vieira by supervising doctoral thesis \cite{S-12}, 
Prof. Peter Wong for pointed the paper \cite{D-G-90} and  Prof. Michael Kelly for help in to improving the exposition this paper.





\bibliographystyle{model1a-num-names}
\bibliography{<your-bib-database>}


\end{document}